\documentclass[11pt]{amsart}
\usepackage{graphicx}
\usepackage{color}
\usepackage{latexsym}
\usepackage{amscd}
\usepackage{a4wide}
\usepackage{epsfig}
\usepackage{float}
\usepackage{amsmath, amssymb, amsthm, amsfonts, amsgen} 
\usepackage[outdir=./]{epstopdf}
\graphicspath{ {./} }

\usepackage{enumitem}
\usepackage{cite}

\usepackage{color}

\renewcommand{\[}{\begin{equation*}}
\renewcommand{\]}{\end{equation*}}
\newcommand{\grad}{\ensuremath{\nabla}}
\newcommand{\p}{\ensuremath{\partial}}

\newcommand{\vp}{\ensuremath{\varphi}}
\newcommand{\al}{\ensuremath{\alpha}}
\newcommand{\be}{\ensuremath{\beta}}
\newcommand{\ga}{\ensuremath{\gamma}}

\newcommand{\la}{\ensuremath{\lambda}}
\renewcommand{\th}{\ensuremath{\theta}}

\renewcommand{\ss}{\ensuremath{\subset}}
\newcommand{\sm}{\ensuremath{\setminus}}
\newcommand{\ra}{\ensuremath{\rightarrow}}

\newcommand{\mb}[1]{\ensuremath{\;\mbox{#1}\;}}
\newcommand{\cd}{\ensuremath{\cdot}}
\newcommand{\ot}{\ensuremath{\otimes}}%Tensorproduct
\newcommand{\ti}{\ensuremath{\times}}%Crossproduct

\newcommand{\Id}{\ensuremath{\mbox{Id}}}

\newcommand{\cof}{\textnormal{cof}\;}
\renewcommand{\div}{\textnormal{div}\;}

\newcommand{\R}{\mathbb{R}}

\newcommand{\N}{\mathbb{N}}

\renewcommand{\d}{\delta}
\newcommand{\A}{\mathcal{A}}

\newtheorem{de}{Definition}[section]
\newtheorem{thm}[de]{Theorem}
\newtheorem{cor}[de]{Corollary}
\newtheorem{lem}[de]{Lemma}

\newtheorem{re}[de]{Remark}
\allowdisplaybreaks
\begin{document}
\title[Uniqueness and regularity in an incompressible variational problem]{A uniqueness criterion and a counterexample to regularity in an incompressible variational problem} % Article title
\author[Marcel Dengler]{M. Dengler}
\address[M. Dengler]{Fliederweg 1, 72189 Vöhringen, Germany, (Corresponding author).}
\email{marci.dengler@web.de}

\author[Jonathan J. Bevan]{J.J. Bevan}
\address[J. Bevan]{Department of Mathematics, University of Surrey, Guildford, GU2 7XH, United Kingdom.}
\email{j.bevan@surrey.ac.uk}
\keywords{Calculus of Variations, Elasticity, Uniqueness, Counterexample to Regularity} 
\subjclass[2020]{49K10, 49K20, 73C15, 73C50}

%\author{M. Dengler and J. J. Bevan}
\maketitle 
\begin{abstract}
\noindent In this paper we consider the problem of minimizing functionals of the form $E(u)=\int_B f(x,\nabla u) \,dx$ in a suitably prepared class of incompressible, planar maps $u: B \to \R^2$.   Here, $B$ is the unit disk and $f(x,\xi)$ is quadratic and convex in $\xi$.  
%certain incompressible variational problems on the full unit ball $B\ss \R^2$ and for a class of quadratic uniformly convex functionals and suitable boundary conditions.  
It is shown that if $u$ is a stationary point of $E$ in a sense that is made clear in the paper, then $u$ is a unique global minimizer of $E(u)$ provided the gradient of the corresponding pressure satisfies a suitable smallness condition. 
%As an application we construct a counterexample to regularity. 
%Indeed, on $B$ and for smooth boundary conditions 
We apply this result to construct a non-autonomous, uniformly convex functional $f(x,\xi)$, depending smoothly on $\xi$ but discontinuously on $x$, whose unique global minimizer is the so-called $N-$covering map, which is Lipschitz but not $C^1$.
\end{abstract}

\section{Introduction} Let $B\ss\R^2$ be the unit ball. For any $u\in W^{1,2}(B,\R^2)$, define the energy $E(u)$ by
\begin{equation}E(u)=\int\limits_B{f(x,\grad u)\;dx},  \label{eq:USPS:1.0}\end{equation}
where the integrand is quadratic in the gradient argument
\begin{align}\label{def:f}f(x,\xi)=M(x)\xi\cd\xi \;\mb{for any} \;x\in B\;\mb{and} \;\xi\in\R^{2\ti2},\end{align}
and where $M\in L^\infty(B,\R^{16})$ is symmetric, i.e. $M_{ijkl}=M_{klij}$ for all $i,j,k,l \in \{1,2\}.  
\footnote{$M\in L^\infty$ is sufficient to guarantee $E(u)<\infty$ for any $u\in W^{1,2}.$}$ Furthermore, we require that there exists a constant $\nu>0$ s.t. 
\begin{equation} M(x)\xi \cdot \xi \geq \nu \, |\xi|^2  \;\mb{for a.e.\@} x\in B \ \textrm{and all }  \xi \in \R^{2\times 2},\label{eq:USPS:1.0.1}\end{equation}
so that $f(x,\xi)$ is uniformly convex in $\xi$.\\

Assume that $g$ is the trace of a map $u_0 \in W^{1,2}(B,\R^2)$ that satisfies $\det \nabla u_0=1$ a.e.\@ in $B$, so that the class \begin{align}\label{eq:1.1.sam}
\A^g:=\{u\in W^{1,2}(B,\R^2): \det\grad u=1 \;\mb{a.e. in}\; B\; u_{|{\p B}}=g\}
\end{align}
is, in particular, nonempty. The constrained minimization problem that we study in this paper is then to find
\begin{equation}
\min\limits_{u\in\A^g} E(u)
\label{eq:USPS:1.1}
\end{equation}
in $\A_g$.  Concrete instances of $g$ for which $\A^g$ is nonempty include:
\begin{itemize}\item[(a)] $g(x):=Ax$, $x\in \partial B$, where $A$ is any constant matrix in $SL^+(2,\R)$, and
\item[(b)] $g(\th):=\frac{1}{\sqrt{N}}e_R(N\th)$,  where $N\in \N\sm\{0\}$ and $0 \leq \th \leq 2\pi$.
\end{itemize}
Note that the latter is the trace of the so-called $N-$covering map 
\begin{align}\label{N-cover}u_{_N}(R,\theta)=\frac{R}{\sqrt{N}}e_{R}(N\theta),\end{align} expressed in plane polar coordinates, and where we employ the notation $e_R(\th)=(\cos \th,\sin\th)$.\\

% and $e_\th(\th)=(-\sin(\th),\cos(\th))$. 

%The conditions set out above allow us, by classical methods, to establish that a minimizer of $E(\cdot)$ in $\A^g$ exists.  
%It is not automatic that a minimizer solves 

%\begin{re}\emph{Note that the above includes the following specific situations:\\
%a) $g(\th)=Ax$ and any constant matrix in $A\in SL^+(2,\R)$ and\\
%b) $ for any $N\in \N\sm\{0\}$ the $N-$covering map on the boundary.}
%\end{re}
\begin{de}(Stationary point)
We say that $u$ is a stationary point of $E(\cdot)$ if there exists a function $\la$, which we shall henceforth refer to as a pressure, belonging to $W^{1,1}(B)$ and such that
\begin{align}\label{def:SP} \div( \nabla_{\xi}f(x,\nabla u) + 2 \la(x) \, \cof \nabla u) = 0 \quad \mathrm{in} \ \mathcal{D}'(B).
\end{align}
\end{de}

The first main result we obtain shows that if $u$ is a stationary point of the energy $E$ whose corresponding pressure $\la$ satisfies, in addition, the assumption that 
\begin{align}\label{ineq:ub}||\grad \la(x)R\|_{L^\infty(B,\R^2)}\le \frac{\sqrt{3}\nu}{2\sqrt{2}},\end{align}
then $u$ is a global minimizer of $E$.  We think of the estimate \eqref{ineq:ub} as characterizing `smallness' of the pressure;  concrete examples (such as can be found in \cite[Proposition 3.4]{JB14} or \cite{BeDe21}, for instance) show that such an estimate need not hold in general.   In the following we assume that $g$ has been chosen and fixed so that $\A^g$ is non-empty.

\begin{thm}[Uniqueness under small pressure]  Let the energy functional $E(u)$ be given by \eqref{eq:USPS:1.0}, and let $f(x,\xi)$ be given by \eqref{def:f}, where $M\in L^{\infty}(B,\R^{16})$ is symmetric and satisfies \eqref{eq:USPS:1.0.1} for some $\nu>0$.  Let $u$ be a stationary point of $E$ in the sense of \eqref{def:SP} and assume that the corresponding pressure $\la$ 
%Let $b\in \R^2,$ $g \in W^{1,2}(S^1,\R^2),$ and let $u\in \A^g$ be an arbitrary stationary point of the energy $E,$ as defined in \eqref{eq:USPS:1.0}.
%Furthermore, assume that for the map $u$ the corresponding pressure $\la$ exists and
satisfies 
\begin{equation}\|\grad \la(x)R\|_{L^\infty(B,\R^2)}\le \frac{\sqrt{3}\nu}{2\sqrt{2}}.\label{eq:Uni.SPC.101}\end{equation}
Then $u$ is a global minimizer of $E$ in $\A^g$.
%problem \eqref{eq:USPS:1.1} is globally minimized by the map $u.$\\

Moreover, if the inequality is strict, i.e. $|\grad \la(x)R|< \frac{\sqrt{3}\nu}{2\sqrt{2}}$ on a measurable set $U\ss B$ with $\mathcal{L}^2(U)>0$ then $u$ is the unique global minimizer of E in $\A^g$.\footnote{Here the norm is defined via $\|f\|_{L^\infty(B,\R^2,\mu)}:=\max\{\|f_1\|_{L^\infty(B,\mu)},\|f_2\|_{L^\infty(B,\mu)}\}$ and $\|f_i\|_{L^\infty(B,\mu)}:=\lim_{p\ra\infty}\left(\int\limits_B{|f_i(x)|^p\;d\mu}\right)^{\frac{1}{p}}.$}
\label{Thm:Uni.SPC.1}
\end{thm}

In our second main result we provide an explicit integrand $f(x,\xi)$ of the form \eqref{def:f} whose corresponding energy functional $E$ is minimized in $\A^g$ by the $N-$covering map $u_{_N}$.  Here, $g$ is the trace of $u_{_N}$ as defined in \eqref{N-cover}.  For its construction we make use of Theorem \ref{Thm:Uni.SPC.1}. 
%As a consequence, our example will be the unique global minimizer of some functional, not just some stationary point. 
A novelty of our approach is that, in order to apply Theorem \ref{Thm:Uni.SPC.1}, we develop a method to compute the corresponding pressure explicitly.  
 %In many other situations it might be enough to guarantee, that the pressure exists in a suitable space, but that's not the case here. Such a counterexample in incompressible elasticity seems to be missing in the literature.\vspace{0.7cm}
 
 \begin{thm}[Counterexample to regularity]\label{cor:4.2.Counter.1} Let $g$ be the trace of the $N-$covering map $u_{_{N}}$, let  $N\in\N\sm\{1\}$ and let $a\in\left(N^2-N,N^2+N\right).$  
Define for $x \in B\sm\{0\}$ and $\xi \in\R^{2 \times 2}$ the function
 \[f(x,\xi)=\nu\left[a(e_{R}^T\xi e_{R})^2+(e_{R}^T\xi e_{\th})^2+a(e_{\th}^T\xi e_{R})^2+(e_{\th}^T\xi e_{\th})^2\right],\]
where $ \nu >0$.  
 %Let $\nu>0,$ $b\in\R^2,$ assume $u_0=g+b$ with $g=\frac{1}{\sqrt{N}}e_{NR}$ and $N\in\N\sm\{0,1\}$ on the boundary, and let the energy be given by \eqref{eq:USPS:1.0}, where the integrand $f$ is of the quadratic form
%\[f(x,\xi)=\nu\left[a(e_{R}^T\xi e_{R})^2+(e_{R}^T\xi e_{\th})^2+a(e_{\th}^T\xi e_{R})^2+(e_{\th}^T\xi e_{\th})^2\right],\]
%for any $x\in B\sm\{0\},$ $\xi\in \R^{2\ti2}$ and some $a\in\left(N^2-N,N^2+N\right).$\vspace{0.7cm}
Then the following statements are true:

\begin{enumerate} [label=(\roman*)]
\item There exists  $M\in L^\infty(B,\R^{16})$ such that
\[f(x,\xi)=\nu M(x)\xi\cd\xi\]
for any $x\in B\sm\{0\},$ $\xi\in \R^{2\ti2}$ and where $M$ is of the explicit form\footnote{Here the multiplication is understood through its action on $\xi\in \R^{2\ti2}$ which is given by
\[((a\ot b)(c\ot d))\xi\cd\xi=(a\ot b)_{ij}\xi_{ij} (c\ot d)_{kl}\xi_{kl}\; \mb{for}\;i,j,k,l\in\{1,2\}.\]
}
\begin{align*}M(x)=&a((e_R\ot e_R)(e_R\ot e_R))+((e_R\ot e_\th)(e_R\ot e_\th))&\\
&+a((e_\th\ot e_R)(e_\th\ot e_R))+((e_\th\ot e_\th)(e_\th\ot e_\th))&\end{align*}
and satisfies 
$M(x) \geq \nu \Id$ for any $x\in B\sm\{0\}$ and any $N\in\N\sm\{1\}.$
\item The maps $x\mapsto M(x)$ and $x\mapsto f(x,\xi),$ for any $\xi\in \R^{2\ti2}\sm\{0\},$ are discontinuous at $0.$
\item The maps $x\mapsto M(x)$ and $x\mapsto f(x,\xi),$ for any $\xi\in \R^{2\ti2}\sm\{0\},$ belong to
\[W^{1,q}\sm W^{1,2} \;\mb{for any}\;1\le q<2\] with the spaces $(B,\R^{16})$ and $(B)$ respectively.
\item The map \begin{equation}u=u_{_N} \in C^{0,1}(B,\R^2)\sm C^1(B,\R^2)\label{eq:Uni.CTR.1.2}\end{equation}
is a stationary point of $E,$ as defined in \eqref{def:SP}, and the corresponding pressure $\la$ exists and satisfies $\la\in W^{1,q}(B) \;\mb{for any}\;1\le q<2.$ 
\item Moreover, the map $u_{_N}$ is the unique global minimizer of $E$ in the class $\A^g.$
\item The minimum energy is given by
\[\min\limits_{v\in\A^g}E(v)=\frac{\nu\pi}{2}(1+a)\left(\frac{1}{N}+N\right).\]
\end{enumerate}
\label{Thm:Counter.1}
\end{thm}

The problem of studying a functional of the form $E(u)$ is of interest not least because the regularity and/or uniqueness of minimizers of such incompressible problems cannot necessarily be determined \emph{a priori}.   Concerning uniqueness in the compressible setting, works including but not limited to  John \cite{J72}, Knops and Stuart \cite{KS84}, Sivaloganathan \cite{S86}, Zhang  \cite{Z91}, and Sivaloganathan and Spector \cite{SS18} provide conditions under which the uniqueness of a global minimizer can be expected. By contrast, a striking example given by Spadaro \cite{S08} clearly demonstrates that global minimizers need not be unique, even under full displacement boundary conditions.   When the domain of integration, or reference configuration, is an annulus, a number of papers, including those of John \cite{J72}, Post and Sivaloganathan \cite{PS97}, Taheri \cite{T09}, and Morris and Taheri \cite{MT17,MT19}, address uniqueness.  With the topology of the annulus at their disposal, a multiplicity of solutions/equilibria can be generated by working with certain homotopy classes.  For example, Morris and Taheri \cite{MT17,MT19} consider functionals of the form $W(x,s,\xi)=F(|x|^2,|s|^2)|\xi|^2/2$, with $F\in C^2$, on the annulus $A$ and admissible maps $\A=W_{\mathrm{id}}^{1,2}(A,\R^2)$, and show that there are countably many solutions, with exactly one for each homotopy class. \\[-3mm]

In the homogeneous, incompressible elasticity setting, Knops and Stuart \cite[Section 6]{KS84} (see also \cite{ST10}) show that the equilibirum solutions they consider are, when subject to affine displacement boundary conditions, global minimizers of the associated energy.   Recent results \cite{BeDe20} show that there are polyconvex energies with inhomogeneous integrands that, under pure displacement boundary conditions, possess countably many pairs of equal-energy stationary points.  It is an open question whether the lowest-energy pair of such stationary points represent global minimizers.    \\[-3mm]

The regularity of equilibrium solutions or mininimzers in incompressible variational problems is, like its compressible counterpart, a delicate matter.   Ball \cite[\S 2.6]{BallOP} points out that any minimizer of the Dirichlet energy in the class $W^{1,2}(B;\R^2)$, subject to the pointwise (incompressibility) constraint $\det \nabla u = 1$ a.e\@ and boundary condition $u(1,\th)=\frac{1}{\sqrt{2}}(\cos(2\th),\sin(2\th))$, must fail to be $C^1$.   There is evidence to suggest that the double-covering map $u_{_2}$ (see \eqref{N-cover}) itself may be the global minimizer in that particular problem  \cite{JB14,BeDe21}.  
A partial regularity result for Lipschitz minimizers that are subject to a type of monotonicity condition was established in \cite{EvGa99}, and Karakhanyan \cite{K11,K13} proves that, in the case of the Dirichlet energy, bounds on the so-called dual pressure lead, by a sophisticated argument, to the conclusion that suitably defined equilibrium solutions must belong to the H\"{o}lder class $C_{\mathrm{loc}}^{\frac{1}{2}}$.  It is perhaps significant that Karakhanyan's results, like ours, also depend on `pressure bounds', but more significant still that the maps he deals with are measure-preserving. The double-covering map $u_2$ mentioned above, and indeed the $N-$covering maps which form the basis of the counterexample to regularity in Section 3 of this paper, do not preserve $\mathcal{L}^2-$measure in the sense of \cite[Eq.\@(24)]{Sv88}, for example, and so are less relevant to physically realistic models of elasticity. \\[-3mm]

%It is common practice to consider the minimization problem \eqref{eq:USPS:1.1} with the set $\A^g$ as given in \eqref{eq:1.1.sam} and the constraint being $\det\grad u=1$  a.e. in $B.$ Notice, however, that in case of $\mb{deg} g\ge 2$ (which is the case for the N-covering map) then $\det\grad u=1$ does no longer represent $\La^2-$ measure preservability. The latter property can be restored by instead of the original minimization problem one instead considers the following: In order for a map $v\in W^{1,2}(B,\R^2)$ with $\mb{tr} (v)=g$ satisfying $\mb{deg} g= N$ for some $N\in N\sm\{0\}$ to be measure preserving one needs that $g$ is the trace of a map $u_0 \in W^{1,2}(B,\R^2)$ that satisfies $\det \nabla u_0=N$ a.e.\@ in $B$. Hence, the set of admissible measure preserving maps can be defined via
%\begin{align*}
%\hat{\mathcal{A}}^{\hat{g}}:=\{v\in W^{1,2}(B,\R^2): \;\det \nabla v=N\; \mb{a.e. in}\; B \; v_{|{\p B}}={\hat{g}}\}.
%\end{align*}
%Finally, realise that for every $u\in\A^g$ the scaled map $v:=\sqrt{N}u\in \hat{\mathcal{A}}^{\hat{g}},$ and vice versa describing a bijection among those maps.  The energy then scales like $E(v)=N\cd E(u)$ implying
%\begin{equation*}
%\min\limits_{v\in\hat{\mathcal{A}}^{\hat{g}}} E(v)=N \min\limits_{u\in\A^g} E(u).
%\label{eq:USPS:1.1}
%\end{equation*}

% Hence the solutions/ minimizers and the corresponding energies of these problems agree up to the given scaling. \\

  It seems that pressure regimes can be used to divide the sorts of incompressible problems we consider into two classes.  The double-covering  problem introduced by Ball appears to lie in the `high pressure' regime\footnote{By which we mean that the pressure $\lambda_{_2}$, say, appearing in the equilibrium equations associated with $u_{_2}$ obeys $||R\,\nabla \la_{_2}||_{\infty}=3\nu$, when adapted to the notation we use in this paper.  The prefactor of $\nu$ in the latter exceeds the prefactor $3^{\frac{1}{2}}2^{-\frac{3}{2}}$ appearing in the condition \eqref{eq:Uni.SPC.101} of Theorem \ref{Thm:Uni.SPC.1}, which is why we refer to this as the `high pressure' regime.},
 whereas the problem we focus on falls, by design, into the `low pressure' regime, where we can say a bit more.   \\[-3mm]
 
Let $v, u \in \mathcal{A}^g$ and suppose that $u$ is a stationary point of $E$ in the sense of \eqref{def:SP}.   To compare $E(v)$ and $E(u)$ we set $\eta=v-u$ and expand $E(v)=E(u+\eta)$ as 
$$E(v)=E(u) + E(\eta) + \int_B 2 M(x) \, \nabla u(x) \cdot \nabla \eta(x) \, dx $$ 
Our problem, as expressed in \eqref{eq:USPS:1.1}, is made more tractable by the observation made in \cite{JB14} that
the stationarity condition \eqref{def:SP} allows us, at the expense of incorporating a pressure term, to rewrite the final, affine-in-$\nabla \eta$ term in the expansion above as a term that is quadratic in $\nabla \eta$, namely
$$ \int_B M(x) \, \nabla u(x)  \cdot \nabla \eta(x) \, dx = \int_B \lambda(x) \det \nabla \eta(x) \, dx.$$
In particular, 
$$E(v) = E(u) + \int_B |\nabla \eta|^2 + 2\lambda \det \nabla \eta \,dx.$$
For the details, see \eqref{eq:Uni.SPC.1} and the foregoing discussion.

\vspace{3mm}

\textbf{Plan of the paper:} The main purpose of Section \ref{sect:uniqueness} will be to prove the uniqueness result, Theorem \ref{Thm:Uni.SPC.1}. We begin by giving two technical lemmas, the first of which enables us to decompose certain expressions in terms of Fourier modes.  Section  \ref{sect:uniqueness} concludes with the proof of Theorem \ref{Thm:Uni.SPC.1}, together with an argument which shows that the prefactor $3^{\frac{1}{2}}2^{-\frac{3}{2}} \simeq 0.6123$ appearing in \eqref{eq:Uni.SPC.101} can be replaced by $1$ when $\la$ depends on just one of the variables $R$, $\theta$.   See Corollary \ref{cor:spudulike}.
%\frac{\sqrt{3}}{2\sqrt{2}},$
%, we prove Theorem \ref{Thm:Uni.SPC.1}. 
%This will be followed by a corollary, where we show that it is possible to improve the prefactor $\frac{\sqrt{3}}{2\sqrt{2}},$ which appears on the RHS of \eqref{eq:Uni.SPC.101}, to 1 if $\la$ only depends on one of its variables $(R,\th).$ 
The focus of Section \ref{Ncoverexample} is Theorem \ref{Thm:Counter.1}.   In order to obtain this result we first develop a method to compute the pressure explicitly:  this is done for a quite general situation in Lemma \ref{lem:Uni.RP.1}, and then more concretely in Lemmata \ref{Lem:Uni.NC.1}-\ref{Lem:Uni.NC.2}.
%, all of which contribute to the proof of Theorem \ref{Thm:Counter.1}.

\subsection{Notation} For a $2\times2-$matrix $A$ the cofactor is given by 
\begin{equation}
\cof A=\begin{pmatrix}a_{22} &-a_{21}\\-a_{12} & a_{11}\end{pmatrix},
\label{eq:1.3}
\end{equation}
and we define the matrix $J$ via \begin{equation*}
J:=\begin{pmatrix}0 &-1\\1 & 0\end{pmatrix}.
\end{equation*}
For two vectors $a\in \R^n,b\in \R^m$ we define the tensor product $a\ot b\in \R^{n\ti m}$ by $(a\ot b)_{ij}:=(ab^T)_{ij}=a_ib_j$ for all $1\le i\le n,$ $1\le j\le m.$  When $\varphi$ is suitably differentiable, we recall that $\det \nabla \varphi = J \varphi_{_R} \cdot \varphi_{_\tau}$, where $\varphi_{_R}$ and $\varphi_{_\tau}=\frac{1}{R}\varphi_{_\th}$ are, respectively, the radial and angular derivatives of $\varphi$.   We use $\mathcal{L}^2$ to denote two-dimensional Lebesgue measure. For any $k\in\N\sm\{0\}$ and $f:B\rightarrow\R$ measurable we define the norm 
$\| f\|_{L^2(dx/R^k)}:=\left(\int\limits_B{|f(x)|^2\;\frac{dx}{R^k}}\right)^{\frac{1}{2}}.$ For a measurable vector-valued $f=(f_1,f_2):B\rightarrow\R^2$ we define $\| f\|_{L^2(dx/R^k)}:=\left(\int\limits_B{|f_1(x)|^2+|f_2(x)|^2\;\frac{dx}{R^k}}\right)^{\frac{1}{2}}.$
All other notation is either standard or is defined when it is used.

\section{Uniqueness in the small pressure regime}\label{sect:uniqueness}

To prove Theorem \ref{Thm:Uni.SPC.1} we need two technical lemmas. The first contains basic properties of functions in the class $W^{1,1}(B)$ that  satisfy $\|R\grad\la\|_{L^\infty(B,\R^2)}<\infty,$ and it relies on a standard Fourier decomposition which, when applied to $\eta\in C^\infty(B,\R^2)$, is given by:
\begin{align*}
\eta(x)=\sum\limits_{j\ge0}\eta^{(j)}(x),\;\mb{where}\;\eta^{(0)}(x)=\frac{1}{2}A_0(R), \;A_0(R)=\frac{1}{2\pi}\int\limits_{0}^{2\pi}{\eta(R,\th)\;d\th}\end{align*} 
and, for any $ j\ge1$,
\begin{align*}
 \eta^{(j)}(x)=A_j(R)\cos(j\th)+B_j(R)\sin(j\th),
 \end{align*}
  where
\begin{align*}
A_j(R)=\frac{1}{2\pi}\int\limits_{0}^{2\pi}{\eta(R,\th)\cos(j\th)\;d\th}\;\mb{and}\;
B_j(R)=\frac{1}{2\pi}\int\limits_{0}^{2\pi}{\eta(R,\th)\sin(j\th)\;d\th}.
\end{align*}
For later use, we set $\tilde{\eta}:=\eta-\eta^{(0)}$.
%\footnote{Note, that this is true for an arbitrary $\eta\in C^\infty(B,\R^2).$ It is not necessary that $\eta$ separates the variables $R$ and $\th.$}\\

\begin{lem}\label{lem:Fdecomp}  Let $\la\in W^{1,1}(B)$ and assume that $\|R\grad\la\|_{L^\infty(B,\R^2)}<\infty.$ 
Then the following statements are true:\\

i) $\la\in BMO(B).$\\
ii) If $\vp_n\ra\vp\in W^{1,2}(B,\R^2)$ then $\int\limits_B{\la(x)\det\grad\vp_n\;dx}\ra\int\limits_B{\la(x)\det\grad\vp\;dx}.$\\
iii) It holds $\int\limits_B{|\grad\vp|^2\;dx}=\sum\limits_{j\ge0}\int\limits_B{|\grad\vp^{(j)}|^2\;dx}$ for any $\vp\in W^{1,2}(B,\R^2).$\\ 
iv) $\det\grad\vp^{(0)}=0$ for any $\vp\in W^{1,2}(B,\R^2).$\\
v)  $\int\limits_B{\la(x)\det\grad\vp\;dx}=-\frac{1}{2}\int\limits_B{((\cof\grad\vp)\grad\la(x))\cd\vp\;dx}$ for any $\vp\in W_0^{1,2}(B,\R^2).$\\
vi)$\int\limits_B{\la(x)\det\grad\vp\;dx}=-\frac{1}{2}\int\limits_B{((\cof\grad\vp^{(0)})\grad\la(x))\cd\tilde{\vp}\;dx}-\frac{1}{2}\int\limits_B{((\cof\grad\vp)\grad\la(x))\cd\tilde{\vp}\;dx}$ for any $\vp\in W_0^{1,2}(B,\R^2).$
\label{Lem:Uni.Tech.1}
\end{lem}
\begin{proof}
\begin{itemize}
\item [i)-iv)] For a proof of these points, see \cite[Lem 3.2 and Prop 3.2]{JB14}. The argument given there still applies if $\la$ depends on $x$ instead of $R.$
\item [v)] Assuming $\vp\in C_c^\infty(B)$, a  computation shows: 
\begin{align*}\int\limits_B{\la(x)\det\grad\vp\;dx}=&\int\limits_B{\la(x)J\vp,_R\cd\vp,_\th\;\frac{dx}{R}}&\\
=&-\int\limits_B{(\la(x)J\vp,_R),_\th\cd\vp\;\frac{dx}{R}}&\\
=&-\int\limits_B{\la(x),_\th J\vp,_R\cd\vp\;\frac{dx}{R}}-\int\limits_B{\la(x)J\vp,_{R\th}\cd\vp\;\frac{dx}{R}}&\\
=&-\int\limits_B{\la(x),_\th J\vp,_R\cd\vp\;\frac{dx}{R}}+\int\limits_B{\la(x),_RJ\vp,_{\th}\cd\vp\;\frac{dx}{R}}&\\
&+\int\limits_B{(\la(x)J\vp,_{\th})\cd\vp,_R\;\frac{dx}{R}}&\\
=&-\int\limits_B{((\cof\grad\vp)\grad\la(x))\cd\vp\;dx}-\int\limits_B{(\la(x)J\vp,_{R})\cd\vp,_\th\;\frac{dx}{R}}&
\end{align*}
The result follows by bringing the rightmost term to the left-hand side and dividing by two.  Note, as a last step, that one needs to upgrade the above equation to hold not just for $\vp\in C_c^\infty(B)$ but instead for all $\vp\in W_0^{1,2}(B).$ This is slightly delicate because of the weak spaces involved: for a proof, see \cite[Lem 3.2.(iv)]{JB14}.\\[-3mm]

\item [vi)] This is a version of (v) in which we emphasise the way that the above expression depends on $\vp^{(0)}$.
Again, we assume $\vp\in C_c^\infty(B),$ and we start by noting $\vp,_\th=\tilde{\vp},_\th,$ hence,
\begin{align*}\int\limits_B{\la(x)\det\grad\vp\;dx}=&\int\limits_B{\la(x)J\vp,_R\cd\tilde{\vp},_\th\;\frac{dx}{R}}&\\
=&-\int\limits_B{((\cof\grad\vp)\grad\la(x))\cd\tilde{\vp}\;dx}+\int\limits_B{(\la(x)J\tilde{\vp},_{\th})\cd\tilde{\vp},_R\;\frac{dx}{R}}&
\end{align*}
then the rightmost term is just
\begin{align*}\int\limits_B{(\la(x)J\tilde{\vp},_{\th})\cd\tilde{\vp},_R\;\frac{dx}{R}}=-\int\limits_B{(\la(x)J\tilde{\vp},_R)\cd\tilde{\vp},_{\th}\;\frac{dx}{R}}=\frac{1}{2}\int\limits_B{((\cof\grad\tilde{\vp})\grad\la(x))\cd\tilde{\vp}\;dx},
\end{align*}
together with the above we get
\begin{align*}\int\limits_B{\la(x)\det\grad\vp\;dx}=&-\frac{1}{2}\int\limits_B{((\cof\grad\vp^{(0)})\grad\la(x))\cd\tilde{\vp}\;dx}-\frac{1}{2}\int\limits_B{((\cof\grad\vp)\grad \la(x))\cd\tilde{\vp}\;dx}.
\end{align*}
\end{itemize}\end{proof}

%\vspace{0.5cm}
The uniqueness condition will be of the form $\|\grad \la(x)R\|_{L^\infty(B,\R^2)}\le C,$ for some constant $C>0$ and where $\la$ will be the corresponding pressure to some stationary point.  A priori, the condition only guarantees the existence of $\grad \la(x)R$ in a suitable space. In the next lemma we show that $\la$ and $\grad\la$ themselves exist in a suitable space, which, in particular, allows one to make use of the technical lemma above.\\

\begin{lem}Let $\mu:B\ra\R$ be a function satisfying 
\[\|\grad \mu(x)R\|_{L^\infty(B,\R^2)}<\infty.\]
Then $\mu\in W^{1,p}(B,\R^2)$ for any $1\le p<2.$\\
\end{lem}
\begin{proof}
The proof is straightforward. Indeed, it holds that
\[\int\limits_B{|\grad \mu|^p\;dx}\le\|\grad \mu R\|_{L^\infty(B,\R^2)}^p\int\limits_B{R^{1-p}\;\frac{dx}{R}},\]
where the latter integrand is integrable for all $1\le p<2.$ \end{proof}

We are now in a position to prove the main statement of this section. 

\vspace{2mm}
\textbf{Proof of Theorem \ref{Thm:Uni.SPC.1}:}\\[-3mm]

Let $u\in\A^g$ be a stationary point with pressure $\la$, let $v\in \A^g$ be arbitrary and set $\eta:=v-u\in W_0^{1,2}(B,\R^2).$\\[-3mm]

We start expanding the energy via 
\[E(v)=E(u)+E(\eta)+H(u,\eta),\]
where
\[H(u,\eta):=2\int\limits_B{M(x)\grad u\cd\grad\eta\;dx}\]
denotes the mixed terms.\\
Expanding the Jacobian of $\eta$ and exploiting the fact that both $u$ and $v$ satisfy $\det\grad u=\det \grad v =1$ a.e. yields
\[\det\grad\eta=-\cof \nabla u \cd \grad \eta \;\mb{a.e.}\]
By the latter identity and the fact that $(u,\la)$ satisfies \eqref{def:SP}, $H$ can be written as
\begin{equation}H(u,\eta)=2\int\limits_B{\la(x)\det \grad \eta\;dx}.\label{eq:Uni.SPC.1}\end{equation}
By Lemma \ref{Lem:Uni.Tech.1}.(vi) we have
\[H(u,\eta)=-\int\limits_B{(\cof\grad\eta^{(0)}\grad\la(x))\cd\tilde{\eta}\;dx}-\int\limits_B{(\cof\grad\eta\grad\la(x))\cd\tilde{\eta}\;dx}=:(I)+(II).\]
Now by noting that the $0-$mode is only a function of $R,$ we get
\begin{align*}
(\cof\grad\eta^{(0)}\grad\la(x))\cd\tilde{\eta}=\frac{\la,_\th}{R}(\eta_{1,R}^{(0)}\tilde{\eta}_2-\eta_{2,R}^{(0)}\tilde{\eta}_1).
\end{align*}
Instead of just $\la,_\th$ on the right hand side of the latter equation we would like to have the full gradient of $\la.$ This can be achieved by using the basic relations $e_\th\cd e_\th=1$ and  $e_R\cd e_\th=0$ to obtain 
\begin{align*}
(\cof\grad\eta^{(0)}\grad\la(x))\cd\tilde{\eta}=(\la,_RR e_R+\la,_\th e_\th)\cd(\eta_{1,R}^{(0)}\tilde{\eta}_2-\eta_{2,R}^{(0)}\tilde{\eta}_1)\frac{e_\th}{R}.
\end{align*}
Arguing similarly for (II), and a short computation shows
\begin{align}
H(u,\eta)=&-\int\limits_B(\la,_RR e_R+\la,_\th e_\th)\cd\left[(\tilde{\eta}_1\tilde{\eta}_{2,\th}-\tilde{\eta}_2\tilde{\eta}_{1,\th})\frac{e_R}{R}\right.&\nonumber\\
&+\left. (\tilde{\eta}_2(\eta_{1,R}^{(0)}+\eta_{1,R})-\tilde{\eta}_1(\eta_{2,R}^{(0)}+\eta_{2,R}))e_\th\right]\;\frac{dx}{R}.&
\label{eq:Uni.SPC.2}\end{align}
By Hölder's inequality we get\footnote{Since we are not only interested in a qualitative but rather a quantitative estimate, we need to specify which norm we pick on $\R^2.$  For this the above Hölder estimate is given more carefully by
\begin{align*}\int\limits_B f\cd g\;d\mu=\int\limits_B f_1 g_1+f_2 g_2\;d\mu\le\int\limits_B |f_1|| g_1|+|f_2|| g_2|\;d\mu\le \max\{\|f_1\|_{L^\infty(B,\mu)},\|f_2\|_{L^\infty(B,\mu)}\}\int\limits_B (| g_1|+| g_2|)\;d\mu\end{align*}
This is the reason why we defined the norm of $f$ via $\|f\|_{L^\infty(B,\R^2,\mu)}:=\max\{\|f_1\|_{L^\infty(B,\mu)},\|f_2\|_{L^\infty(B,\mu)}\}.$}
\begin{align*} 
H(u,\eta)\ge&-\|\grad \la(x)R\|_{L^\infty(B,\R^2,\frac{dx}{R})}\int\limits_B\left[\left|\tilde{\eta}_1\tilde{\eta}_{2,\th}-\tilde{\eta}_2\tilde{\eta}_{1,\th}\right|\frac{1}{R}\right.&\\
&+\left.\left|\tilde{\eta}_2(\eta_{1,R}^{(0)}+\eta_{1,R})-\tilde{\eta}_1(\eta_{2,R}^{(0)}+\eta_{2,R})\right|\right]\;\frac{dx}{R}.&
\end{align*}
By $\|\grad\la(x) R\|_{L^\infty(\frac{dx}{R})}\le\frac{\sqrt{3}\nu}{2\sqrt{2}}$ and a weighted Cauchy-Schwarz Inequality, we see
\begin{align*}
H(u,\eta)\ge&-\frac{\nu\sqrt{3}}{4\sqrt{2}}\left[2a\|\tilde{\eta}_1\|_{L^2(dx/R^2)}^2+2a\|\tilde{\eta}_2\|_{L^2(dx/R^2)}^2\right.&\\
&\left.+\frac{1}{a}\int\limits_B{\left[\frac{\tilde{\eta}_{2,\th}^2}{R^2}+(\eta_{2,R}^{(0)}+\eta_{2,R})^2+(\eta_{1,R}^{(0)}+\eta_{1,R})^2+\frac{\tilde{\eta}_{1,\th}^2}{R^2}\right]\;dx}\right].&
\end{align*}
Next we recall an elementary Fourier estimate (see, for instance, \cite[Proof of Proposition 3.3]{JB14}), which states that for any $\phi\in C^\infty(B)$ it holds
\begin{equation}\int\limits_B{R^{-2}|\tilde{\phi}_{,\th}|^2\;dx}\ge\int\limits_B{R^{-2}|\tilde{\phi}|^2\;dx}.\label{eq:Uni.Buckling.1}\end{equation}

Applying the Cauchy-Schwarz inequality and \eqref{eq:Uni.Buckling.1}, and then  combining some of the norms yields
\begin{align*}
H(u,\eta)\ge&-\frac{\nu\sqrt{3}}{4\sqrt{2}}\left[(2a+\frac{1}{a})\|\tilde{\eta}_1,_\th\|_{L^2(dx/R^2)}^2+(2a+\frac{1}{a})\|\tilde{\eta}_2,_\th\|_{L^2(dx/R^2)}^2 +\right.&\\
&\left.+\frac{2}{a}\|\eta_1,_R^{(0)}\|_{L^2(dx)}^2+\frac{2}{a}\|\eta_1,_R\|_{L^2(dx)}^2+\frac{2}{a}\|\eta_2,_R^{(0)}\|_{L^2(dx)}^2+\frac{2}{a}\|\eta_2,_R\|_{L^2(dx)}^2\right]&\\
\ge&-\frac{\nu\sqrt{3}}{4\sqrt{2}}\left[(2a+\frac{1}{a})\|\tilde{\eta},_\th\|_{L^2(dx/R^2)}^2+\frac{2}{a}\|\eta,_R^{(0)}\|_{L^2(dx)}^2+\frac{2}{a}\|\eta,_R\|_{L^2(dx)}^2\right].&
\end{align*}
Making use of $\tilde{\eta},_\th=\eta,_\th,$ which is true since the zero-mode does not depend on $\th,$ and $\|\eta,_R^{(0)}\|_{L^2(dx)}^2\le\|\eta,_R\|_{L^2(dx)}^2$ we obtain
\begin{align*}
H(u,\eta)\ge&-\frac{\nu\sqrt{3}}{4\sqrt{2}}\left[(2a+\frac{1}{a})\|\eta,_\th\|_{L^2(dx/R^2)}^2+\frac{4}{a}\|\eta,_R\|_{L^2(dx)}^2\right].&
\end{align*}
Choosing $a=\frac{\sqrt{3}}{\sqrt{2}}$ and again combining norms gives
\begin{align*}
H(u,\eta)\ge&-\frac{\nu\sqrt{3}}{4\sqrt{2}}\left[\frac{4\sqrt{2}}{\sqrt{3}}(\|\eta,_\th\|_{L^2(dx/R^2)}^2+\|\eta,_R\|_{L^2(dx)}^2)\right]&\\
=&-\nu D(\eta),&
\end{align*}
where $D(\eta):=\|\grad\eta\|_{L^2(dx)}^2$ denotes the Dirichlet energy.
This yields
\[E(\eta)+H(u,\eta)\ge E(\eta)-\nu D(\eta)\ge0,\]
which, since $M\xi \cdot \xi \geq \nu|\xi|^2$ for all $\xi \in \R^2$, completes the proof.

The prefactor $\frac{\sqrt{3}}{2\sqrt{2}}$ in \eqref{eq:Uni.SPC.101} is the best we have for general $\la$ at the moment.  If, however, circumstances are such that $\la$ depends on only one of $R,\th$ throughout $B$, then condition \eqref{eq:Uni.SPC.101} can be replaced by the weaker assumption \begin{align}\label{ineq:ubbetter}\|\grad \la(x)R\|_{L^\infty(B,\R^2)}\le \nu.\end{align}

\begin{cor}\label{cor:spudulike} Let the conditions of Theorem \ref{Thm:Uni.SPC.1} be in force, but with \eqref{ineq:ubbetter} replacing \eqref{ineq:ub}, and assume  
that either $\la(x)=\la(R)$ or $\la(x)=\la(\th)$ for all $x\in B$.   Then the conclusions of Theorem \ref{Thm:Uni.SPC.1} hold.
\end{cor}
\begin{proof}
\textbf{(i)} ($\la(x)=\la(R)$.) 
This case is significantly simpler and one can argue more along the lines of the proof of \cite[Prop.3.3]{JB14}.
The reason is that in this case it still holds that
\[H(v,\eta)=2\int\limits_B{\la(R)\det \grad \eta\;dx}=2\int\limits_B{\la(R)\det \grad \tilde{\eta}\;dx},\]
where $\tilde{\eta}=\eta-\eta^{(0)}$ eliminating the $0-$mode.\\ Then applying  of \cite[Lemma 3.2.(iv)]{JB14} yields
\[H(v,\eta)=\int\limits_B{\la'(R)R\tilde{\eta}\cd J \tilde{\eta}_{,\th}\;\frac{dx}{R^{2}}}.\]
Using Hölder's inequality, $\|\la'(R)R\|_{L^\infty(0,1)}\le\nu,$ and Fourier estimate \eqref{eq:Uni.Buckling.1} we get
\begin{align*}
H(v,\eta)\ge&-\|\la'(R)R\|_{L^\infty(0,1)} \int\limits_B{|\tilde{\eta}||\tilde{\eta}_{,\th}|\;\frac{dx}{R^{2}}}&\\
\ge&-\nu \left(\int\limits_B{|\tilde{\eta}|^2\;\frac{dx}{R^{2}}}\right)^{1/2}\left(\int\limits_B{|\tilde{\eta}_{,\th}|^2\;\frac{dx}{R^{2}}}\right)^{1/2}&\\
\ge&-\nu  \int\limits_B{|\tilde{\eta}_{,\th}|^2\;\frac{dx}{R^{2}}}&\\
\ge&-\nu  \int\limits_B{|\grad \eta|^2\;dx}.&
\end{align*}
Note, as before, that the $\sim$ could be dropped because $ \int\limits_B{|\grad \tilde{\eta}|^2\;dx}\le \int\limits_B{|\grad \eta|^2\;dx}.$ \vspace{0.5cm}

\textbf{(ii)} ($\la(x)=\la(\th).$)
Here we start with \eqref{eq:Uni.SPC.2} which simplifies to
\begin{align*}
H(u,\eta)=&-\int\limits_B{\la,_\th(\th)[(\eta_2,_R^{(0)}+\eta_2,_R)\tilde{\eta}_1+(\eta_1,_R^{(0)}+\eta_1,_R)\tilde{\eta}_2]\;\frac{dx}{R}}&
\end{align*}
By Hölder's inequality, Inequality \eqref{eq:Uni.Buckling.1} and  $\|\la,_\th\|_{L^\infty(0,2\pi)}\le\nu$ we get
\begin{align*}
H(u,\eta)\ge&-\|\la,_\th\|_{L^\infty(0,2\pi)}\int\limits_B{|(\eta_2,_R^{(0)}+\eta_2,_R)\tilde{\eta}_1+(\eta_1,_R^{(0)}+\eta_1,_R)\tilde{\eta}_2|\;\frac{dx}{R}}&\\
\ge&-\frac{\nu}{2}[2\|\tilde{\eta}_1,_\th\|_{L^2(dx/R^2)}^2+2\|\tilde{\eta}_2,_\th\|_{L^2(dx/R^2)}^2+\|\eta,_R^{(0)}\|_{L^2(dx)}^2+\|\eta,_R\|_{L^2(dx)}^2]&
\end{align*}
Using $\tilde{\eta},_\th=\eta,_\th$ and $\|\eta,_R^{(0)}\|_{L^2(dx)}^2\le\|\eta,_R\|_{L^2(dx)}^2$ we get
\begin{align*}
H(u,\eta)\ge&-\frac{\nu}{2}[2\|\eta,_\th\|_{L^2(dx/R^2)}^2+2\|\eta,_R\|_{L^2(dx)}^2]&\\
=&-\nu D(\eta).&
\end{align*}
\end{proof}

\begin{re}[Relaxation of the assumptions]
The result of Theorem \ref{Thm:Uni.SPC.1} continues to hold if we assume that $f(x;\xi)=M(x)\xi \cdot \xi \geq \nu(|x|)|\xi|^2$ for some $\nu\in L^\infty(\R_{+})$, $\nu(R) \geq 0$, and all $\xi \in \R^{2\times 2}$. 
Here, $\nu(R)=0$ is allowed\footnote{Note that in these circumstances, we do not need to verify that a minimizer of the associated functional $E(u)$ exists in order to apply Theorem \ref{Thm:Uni.SPC.1}.  Rather, it is enough to establish that $u$ is a stationary point in the sense of \eqref{def:SP}.}  to be $0$. The assumption $\nu=\nu(R)$ is needed because we do not know if the Fourier estimate \eqref{eq:Uni.Buckling.1} is still true if $\nu$ depends on both $R$ and $\th.$ 
In this case, the small pressure condition can be relaxed to a pointwise estimate:
\[|\grad \la(x)R|\le \frac{\sqrt{3}\nu(R)}{2\sqrt{2}} \quad \mathrm{for \ a.e.} \ x \in B,\]
with uniqueness if the inequality is strict on some non-null set.
\end{re}

\section{A method for computing the pressure and a counterexample to regularity}\label{Ncoverexample}
In this section we a construct an explicit functional $E(u)$ of the form $\eqref{eq:USPS:1.0}$, where the integrand obeys \eqref{def:f} and \eqref{eq:USPS:1.0.1}, such that $u=u_{_N}$ is the global minimizer of $E$ in $\mathcal{A}^{\mathrm{tr}\,(u_{_{\small{N}}})}$.   The strategy is as follows:  
\begin{itemize}\item[(i)] select a candidate trace function $g \in W^{k+1,p}(\mathbb{S}^1,\R^2)$ for $k \geq 1$ and $1 \leq p \leq \infty$ which obeys\footnote{The condition $Jg\cdot g'=1$ ensures that the one-homogeneous extension $u$ obeys $\det \nabla u =1$ a.e.\@ in $B$} $1=Jg(\theta)\cdot g'(\theta)$ for a.e. $\theta$ in $[0,2\pi)$\footnote{Here, by a slight abuse of notation, we put $\tilde{g}(\theta)=g(\cos\th,\sin \th)$ and then promptly drop the \ $\tilde{\!}$}
\item[(ii)] extend $g$ to a one-homogeneous function $u(R,\th):=R g(\theta)$, and compute, in Lemma \ref{lem:Uni.RP.1}, a PDE which must be satisfied by both $u$ and an associated $\lambda$ in order that $u$ is a stationary point of $E$ in the sense of \eqref{def:SP};
\item[(iii)] fix $g=\mathrm{tr}\,(u_{_{\small{N}}})$ and construct, in Lemma \ref{Lem:Uni.NC.1}, a suitable $f(x,\xi)$ such that the PDE in step (ii) can be solved for $\lambda$, and 
\item[(iv)] verify that the small pressure condition stipulated in Corollary \ref{cor:spudulike} is satisfied by $\lambda$, and hence that $u_{_{\small{N}}}$ is the unique global minimizer of the associated energy $E$. 
\end{itemize}

We first examine conditions on $M$ necessary for a to be a stationary point of $E$ in the sense of \eqref{def:SP}.  

%How does this relate to the small pressure criterion discussed in the previous section? In order to understand the idea, realise firstly that $g=\frac{1}{\sqrt{N}}e_{R}(N\th)$ is a possible choice as boundary conditions and secondly that $u=\frac{R}{\sqrt{N}}e_{R}(N\th)+b\in C^{0,1}(B,\R^2)\sm C^{1}(B,\R^2)$ for any $N\in\N\sm\{0,1\}$ and any $b\in \R^2.$ Now if we are able to construct a concrete $M$ s.t. the corresponding pressure $\la$ satisfies the small pressure criteria, then we could guarantee that the above map has to be a global minimizer of the corresponding energy $E.$\\
%For this sake, we start by considering 1-homogenous stationary points of the functional $E,$ as defined in \eqref{eq:USPS:1.0}, and we compute the corresponding quantity $\grad \la(x)R.$ The next lemma shows that $\grad \la(x)R$ needs to satisfy a PDE system weakly and that $\grad \la(x)R$ is completely determined by $M$ (corresponding to the functional) and the boundary conditions $g.$\\

\textbf{Notation:} Recall the notation for 2d polar coordinates \[\{e_R,e_\th\}:=\{(\cos \th,\sin\th),(-\sin\th,\cos\th)\}.\] Additionally, we will use \[\{e_{NR},e_{N\th}\}:=\{(\cos(N\th),\sin(N\th)),(-\sin(N\th),\cos(N\th))\}\] for any $N\in \N.$ Moreover, we will use the notation $M_{ijgk}=(M(e_i\ot e_j))\cd(g\ot e_k)$ for any combination of $i,j,k\in\{R,\th\}$ and any map $g \in \R^{2}.$ Especially, if $g=e_{Nl}$ for some $l\in\{R,\th\}$ we will use $M_{ij(Nl)k}$ for short.\\[-3mm]

\begin{lem}[Representation of the pressure] Let $1\le p\le\infty,$ $k\in\N\sm\{0\}$ and assume $M\in (L^{\infty}\cap W^{k,p})(B,\R^{16}),$ $g \in W^{k+1,p}(\mathbb{S}^1,\R^2)$  where $g$ obeys $Jg \cdot g'=1$ a.e. in $[0,2\pi),$ and let $u=Rg(\th) \in \A^g$ be a stationary point of the energy $E$ as defined in \eqref{def:SP}.\\

Then there exists a corresponding pressure $\la\in W^{k,p}(B,\R)$ and it satisfies the following system of equations a.e.\@ in $B:$
\begin{align}
\la(x),_\th(Jg\cd e_R)-\la(x),_RR(Jg'\cd e_R)=&-[M_{R\th (g+g'')\th}&\nonumber\\
&+((M,_\th)_{R\th g R}+(M,_\th)_{R \th g'\th})&\nonumber\\
&+ R((M,_R)_{R R g R}+(M,_R)_{R R g' \th})]&\nonumber\\
=:&h_1(M,g)&\nonumber\\
\la(x),_\th(Jg\cd e_\th)-\la(x),_RR(Jg'\cd e_\th)=&-[M_{\th \th (g+g'')\th}&\nonumber\\
&+((M,_\th)_{\th \th g R}+(M,_\th)_{\th \th g'\th})&\nonumber\\
&+ R((M,_R)_{\th R g R}+(M,_R)_{\th R g' \th})]&\nonumber\\
=:&h_2(M,g)&
\label{eq:Uni.RP.1}
\end{align}
\label{lem:Uni.RP.1}
\end{lem}
\begin{proof}
Let $u=Rg(\th)\in \A^g$ be a stationary point. If there exists a corresponding pressure $\la\in W^{1,p}$ then $u$ is a solution of
\begin{equation}\int\limits_B{M(x)\grad u\cd\grad\eta\;dx}=-\int\limits_B{\la(x)\cof\grad u\cd\grad\eta\;dx}\;\mb{for any} \eta\in C_c^\infty(B,\R^2).\label{eq:Uni.RP.2}\end{equation}

For now, let us assume that $\la\in W^{1,p}(B,\R).$ In order to derive the system of equations above, we enter the explicit form of $u$ and the representation $\eta=(\eta\cd e_R) e_R+(\eta\cd e_{\th}) e_{\th}$ into the stationarity condition. By some further calculations, which are mainly integrations by parts, we obtain \eqref{eq:Uni.RP.1}.  In the last step of the proof we discuss the existence of $\la\in W^{1,p}.$\\

\textbf{Step 1:} Computation of left-hand side of \eqref{eq:Uni.RP.2}:\\
The derivative and the cofactor of the map $u=Rg(\th)+b$ are given by
\[\grad u=g\ot e_R+g'\ot e_\th,\]
\[\cof\grad u=Jg\ot e_\th-Jg'\ot e_R.\]

Plugging the above into the left-hand side of \eqref{eq:Uni.RP.2} and integrating by parts yields
\begin{align*}(LHS)=&\int\limits_B{M(x)(g\ot e_R+g'\ot e_\th)\cd\left(\eta,_R\ot e_R+\frac{1}{R}\eta,_\th\ot e_\th\right)\;dx}&\\
=&-\int\limits_B{R(M(x),_R)(g\ot e_R+g'\ot e_\th)\cd(\eta\ot e_R)\;\frac{dx}{R}}&\\
&-\int\limits_B{[M(x)((g+g'')\ot e_\th)+M(x),_\th(g\ot e_R+g'\ot e_\th)]\cd(\eta\ot e_\th)\;\frac{dx}{R}}&
\end{align*}
Now by expanding  $\eta=\al e_R+\be e_{\th}$ with $\al=(\eta\cd e_R)$ and $\be=(\eta\cd e_{\th})$ and the shorthand introduced above we get
\begin{align*}(LHS)=&-\int\limits_B{R[(M,_R)_{RRgR}+(M,_R)_{RRg'\th}]\al+R[(M,_R)_{\th RgR}+(M,_R)_{\th Rg'\th}]\be\;\frac{dx}{R}}&\\
&-\int\limits_B{M_{R\th(g+g'')\th}\al+[(M,_\th)_{R\th gR}+(M,_\th)_{R\th g'\th}]\al}&\\
&{+M_{\th\th(g+g'')\th}\be+[(M,_\th)_{\th\th gR}+(M,_\th)_{\th\th g'\th}]\be\;\frac{dx}{R}}&\\
=&\int\limits_B{h_1\al+h_2\be\;\frac{dx}{R}}.&
\end{align*}

\textbf{Step 2:} Computation of right-hand side of \eqref{eq:Uni.RP.2}:  \\
Now by again using the explicit form of $\cof\grad u,$ and integration by parts we get
\begin{align*}
(RHS)=&-\int\limits_B{(\la(x)(Jg\cd\eta,_\th)-\la(x)R(Jg'\cd\eta,_R))\;\frac{dx}{R}}&\\
=&\int\limits_B{\la,_\th(x)(Jg\cd\eta)-\la,_R(x)R(Jg'\cd\eta)\;\frac{dx}{R}}.&
\end{align*}
Further, we use the expression $\eta=\al e_R+\be e_{\th}$ with the notation $\al=(\eta\cd e_R)$ and $\be=(\eta\cd e_{\th})$ to derive
\begin{align*}
(RHS)=&\int\limits_B{(\la,_\th(x)(Jg\cd e_R)-\la,_R(x)R(Jg'\cd e_R))\al}&\\
&{+(\la,_\th(x)(Jg\cd e_\th)-\la,_R(x)R(Jg'\cd e_\th))\be\;\frac{dx}{R}}.&
\end{align*}
Together with Step 1 and the realization that in the above $\al,\be\in C_c^\infty(B)$ are arbitrary, the claimed equations need to be true a.e. in $B$.\\[-3mm]

\textbf{Step 3:} Existence of the pressure $\la\in W^{1,p}(B,\R):$\\
The equations above can be rewritten as 
\[\int\limits_B{\div(\la \, \cof\grad u)\cd\eta\;dx}=\int\limits_B{h(M,g)\cd\eta\;\frac{dx}{R}},\]
where $h=(h_1,h_2).$ We know that $\div(\la\cof\grad u)\in L^p(dx)$ iff $h(M,g)\in L^p(\frac{dx}{R}),$ with obvious notation.  Now consider $h_1(M,g)$ (the argument being similar for $h_2$) and define
\begin{align*}
h_{11}:=&-M_{R\th (g+g'')\th}&\\
h_{12}:=&[((M,_\th)_{R\th g R}+(M,_\th)_{R \th g'\th})&\\
&+ R((M,_R)_{R R g R}+(M,_R)_{R R g' \th})].&
\end{align*}
Then for $h_{11}\in L^p(\frac{dx}{R})$ we need $M\in L^\infty(dx)$ and $g,g''\in L^p,$ which is true by assumption. Now, by Sobolev imbedding, we have $W^{2,p}\hookrightarrow W^{1,\infty}([0,2\pi),\R^2),$ and hence, in order for $h_{12}\in L^p(\frac{dx}{R}),$ and bearing in mind that $g,g'\in L^\infty,$ it is enough to require that $\grad M\in L^{p}(dx).$
This is exactly how we chose the classes for $M$ and $g.$  This guarantees the existence of $\div(\la\cof\grad u)=(\cof\grad u)\grad\la\in L^p(dx).$ By further noting that $g,g'\in L^\infty$, it is immediate that $\grad u\in L^\infty(dx),$ and since $\det \nabla u=1$ a.e.\@ in $B$, we may write 
$$ \grad\la=(\cof \nabla u)^T \frac{h(M,g)}{R} \in L^p(dx).$$
%\grad u^T(\cof\grad u)\grad\la\in L^p(dx)$$
%Now by the relation $A^T\cof A=(\det A)\Id,$ which is true for any $A\in \R^{2\ti2},$ we have $\grad\la=\grad u^T(\cof\grad u)\grad\la\in L^p(dx).$ 
In particular, when $M$ and $g$ are specified, $\nabla \lambda$ is specified and it belongs to the class $L^p(B)$, reverting to the traditional notation.  One can argue similarly for the higher integrability. \end{proof}

We now specify $g=\mathrm{tr}\,(u_{_{\small{N}}})$ and compute the pressure under the assumptions that $M$ depends only on $\th$ and is diagonal with respect to the basis of polar coordinates.\\[-3mm]
 
\begin{lem}[Representation of the pressure, N-cover, M(\th)=diag] For $N\in\N\sm\{1\}$ let $g=\frac{1}{\sqrt{N}}e_{NR}$ and assume $M\in (L^{\infty}\cap W^{k,p})(B,\R^{16})$ for some $1\le p\le\infty$ and $k\in\N,$ where $M$ is of the specific form
\[M(x)=diag(M_{R R R R},M_{R\th R \th},M_{\th R \th R},M_{\th \th \th \th})=diag(\al(\th),\be(\th),\ga(\th),\d(\th))\] with $\nu>0$ and $\al,\be,\ga,\d\ge\nu$ for any $\th\in[0,2\pi).$ Furthermore, suppose $u=Rg(\th)\in \A^g$ is a stationary point of the energy $E,$ as defined in \eqref{eq:USPS:1.0}.\\
Then there exists a corresponding pressure $\la\in W^{k,p}(B)$ and it satisfies the following system of equations a.e.\@ in $B:$
\begin{align}
-\la(x),_\th\frac{1}{\sqrt{N}}\sin(\th_{N-1})+\la(x),_RR\sqrt{N}\cos(\th_{N-1})=&\sqrt{N}\be'\sin(\th_{N-1})&\nonumber\\
+\left[\sqrt{N}(N-1)\be+\sqrt{N}\d-\right.&\left.\frac{\al}{\sqrt{N}}\right] \cos(\th_{N-1})&\nonumber\\
=:&h_1& \label{eq:Uni.NC.1}\\
\la(x),_\th\frac{1}{\sqrt{N}}\cos(\th_{N-1})+\la(x),_RR\sqrt{N}\sin(\th_{N-1})=&-\sqrt{N}\d'\cos(\th_{N-1})&\nonumber\\
+\left[\sqrt{N}\be+\sqrt{N}(N-1)\d-\right.&\left.\frac{\ga}{\sqrt{N}} \right] \sin(\th_{N-1})&\nonumber\\
=:&h_2&
\label{eq:Uni.NC.2}
\end{align}
where we used the shorthand $\th_{k}:=k\th$ for any $k\in\R.$ 
\label{Lem:Uni.NC.1}
\end{lem}
\begin{proof}
By Lemma \ref{lem:Uni.RP.1}, we know that the pressure $\la$ exists and system \eqref{eq:Uni.RP.1} is satisfied. Now we just have to verify that \eqref{eq:Uni.RP.1} agrees with the claimed system given by \eqref{eq:Uni.NC.1}\ and \eqref{eq:Uni.NC.2}. We start by verifying the first of the equations in the system \eqref{eq:Uni.RP.1}.\\

\textbf{Step 1:}\\
That the left-hand side of \eqref{eq:Uni.NC.1} follows from the left-hand side of the first equation in system \eqref{eq:Uni.RP.1} is a straightforward calculation. Hence, we focus on the corresponding right-hand side, which we named $h_1.$  First note that, because $M$ depends only on $\th$, we are left with
\begin{equation}h_1=-\left[M_{R\th (g+g'')\th}+(M,_\th)_{R\th g R}+(M,_\th)_{R \th g'\th}\right].\label{eq:Uni.NC.3}\end{equation}
We have $g=\frac{1}{\sqrt{N}}e_{NR}, g'=\sqrt{N}e_{N\th},g''=-N\sqrt{N}e_{NR}$, and hence
\[M_{R\th (g+g'')\th}=\left(\frac{1}{\sqrt{N}}-\sqrt{N}N\right)M_{R\th (NR)\th}=\left(\frac{1}{\sqrt{N}}-\sqrt{N}N\right)(M_{R\th R\th}(e_{NR}\cd e_R)+M_{R\th \th \th}(e_{NR}\cd e_\th)).\]
Using that $M_{R\th \th \th}=0$ and $M_{R\th R\th}=\be$ yields
\[M_{R\th (g+g'')\th}=\left(\frac{1}{\sqrt{N}}-\sqrt{N}N\right)\be\cos(\th_{N-1}).\]
For the second term of \eqref{eq:Uni.NC.3}, consider
\[M_{R\th g R,\th}=(M,_\th)_{R\th g R}+M_{\th \th g R}-M_{R R g R}+M_{R\th g' R}+M_{R\th g \th},\]
which, after a short calculation, gives
\begin{align*}
(M,_\th)_{R\th g R}=&M_{R\th g R,\th}-M_{\th \th g R}+M_{R R g R}-M_{R\th g' R}-M_{R\th g \th}&\\
=&\frac{1}{\sqrt{N}}[\al-\be]\cos(\th_{N-1}).&
\end{align*}
Similarly, for the rightmost term of \eqref{eq:Uni.NC.3} we get
\begin{align*}
(M,_\th)_{R\th g' \th}=&M_{R\th g' \th,\th}-M_{\th \th g' \th}+M_{R R g' \th}-M_{R\th g'' \th}+M_{R\th g' R}&\\
=&-\sqrt{N}\be'\sin(\th_{N-1})+\sqrt{N}[\be-\d]\cos(\th_{N-1}).&
\end{align*}
Together,
\[h_1=\sqrt{N}\be'\sin(\th_{N-1})+\left[\sqrt{N}\d-\frac{\al}{\sqrt{N}}+\sqrt{N}(N-1)\be\right]\cos(\th_{N-1}).\]

\textbf{Step 2:}\\
By arguing similarly, we find that
\[h_2=-\left[M_{\th \th (g+g'')\th}+(M,_\th)_{\th \th g R}+(M,_\th)_{\th \th g'\th}\right].\]
Then 
\begin{align*}
M_{\th \th (g+g'')\th}=&\left(\frac{1}{\sqrt{N}}-\sqrt{N}N\right)\d\sin(\th_{N-1})&\\
(M,_\th)_{\th \th g R}=&M_{\th \th g R,\th}+M_{R \th g R}+M_{\th R g R}-M_{\th \th g' R}-M_{\th \th g \th}&\\
=&\frac{1}{\sqrt{N}}[\ga-\d]\sin(\th_{N-1})&\\
(M,_\th)_{\th \th g' \th}=&M_{\th \th g' \th,\th}+M_{R \th g' \th}+M_{\th R g' \th}-M_{\th \th g'' \th}+M_{\th \th g' R}&\\
=&\sqrt{N}\d'\cos(\th_{N-1})+\sqrt{N}[\d-\be]\sin(\th_{N-1})&
\end{align*}
and finally
\[h_2=-\sqrt{N}\d'\cos(\th_{N-1})+\left[\sqrt{N}\be+\sqrt{N}(N-1)\d-\frac{\ga}{\sqrt{N}}\right]\sin(\th_{N-1}),\]
completing the proof.\end{proof}

Next, we compute the small pressure criteria in the same situation. Moreover, we will provide an explicit form of the pressure.

\begin{lem}[Small pressure condition, N-cover, M(\th)=diag] Let the assumptions be as above.
For any $N\in\N\sm\{1\},$ let $M=(a,1,a,1)\nu,$ where we pick $a$ to be constant and in the range
\[1\le N^2-N< a < N^2+N.\]
Then for this $M$ the corresponding pressure $\la$ is given by
\begin{align*}
\la(x)=c+\left[N-\frac{a}{N}\right]\ln(R) \;\mb{for any}\; x\in B
\end{align*}
for any real constant $c\in \R,$ which is independent of $R$ and $\th.$ Moreover, 
$\la\in W^{1,q}(B) \;\mb{for any}\;1\le q<2$ and $\la$ satisfies condition \eqref{eq:Uni.SPC.101} strictly.
\label{Lem:Uni.NC.2}
\end{lem}
\begin{proof}
Define first
\[H_1=\left[\sqrt{N}(N-1)\be+\sqrt{N}\d-\frac{\al}{\sqrt{N}}\right]\;\mb{and}\;H_2=\left[\sqrt{N}\be+\sqrt{N}(N-1)\d-\frac{\ga}{\sqrt{N}}\right].\]
By solving the system \eqref{eq:Uni.NC.1} and \eqref{eq:Uni.NC.2} we obtain
\begin{align*}
\la,_RR=&(\be'-\d')\frac{\sin(2\th_{N-1})}{2}+\frac{1}{\sqrt{N}}(H_1\cos^2(\th_{N-1})+H_2\sin^2(\th_{N-1}))&\\
\la,_\th=&\sqrt{N}(H_2-H_1)\frac{\sin(2\th_{N-1})}{2}-N(\be'\sin^2(\th_{N-1})+\d'\cos^2(\th_{N-1})).&
\end{align*}
For the specific case of $M=(a,1,a,1)\nu$ they become
\begin{align*}
\la,_RR=\left[N-\frac{a}{N}\right] \;\mb{and}\;\la,_\th=0
\end{align*}
showing, in particular, that $\la$ depends only on $R$ i.e. $\la(x)=\la(R).$
Indeed, the pressure is then given by
\begin{align*}
\la(x)=c+\left[N-\frac{a}{N}\right]\ln(R) \;\mb{for any}\; x\in B
\end{align*}
and for any real constant $c\in \R,$ which is independent of $R$ and $\th.$
The small pressure condition of Corollary \ref{cor:spudulike} can now be applied, giving
\begin{align*}
\left|N-\frac{a}{N}\right|<1.
\end{align*}
Solving this inequality by case distinction yields the claimed bounds on $a.$ The integrability is then easily deduced, completing the proof.\end{proof}

To make it more accessible for the reader we collect what we have shown so far in the following.\\

\textbf{Proof of Theorem \ref{Thm:Counter.1}}
\begin{enumerate} [label=(\roman*)]
\item\!\!, (ii), and (vi) trivial.
\item[(iii)]It is enough to show this point for $M.$ Note, that $M$ only depends on $\th,$ i.e. $M(x)=M(\th).$ Hence, the gradient is given by
\begin{equation}\grad M=\frac{1}{R}\p_\th M(\th)\ot e_\th \;\mb{for any}\; x\in B\sm\{0\}.\label{eq:Uni.CTR.1.3}\end{equation}
First realise that the derivative with respect to $\th$ only replaces $e_R$ with $e_\th$ (up to sign) and vice versa, and therefore one can still bound the modulus of $\|\p_\th M(\th)\|_{L^\infty(B,\R^{16})}\le C(a)$ via some real constant $C(a)>0$. Then, by integrating $|\grad M|^q $ with respect to $dx$ using \eqref{eq:Uni.CTR.1.3} and by the latter discussion, the claim follows.

\item[(iv)] As a consequence of $g\in C^\infty$ and point (iii), Lemma \ref{Lem:Uni.NC.1} guarantees that $u$ is a stationary point and the existence of $\la$ in the right spaces.

\item[(v)] By Lemma \ref{Lem:Uni.NC.2} we know that $\la$ satisfies the small pressure criteria strictly. Together, with Theorem \ref{Thm:Uni.SPC.1} this implies that $u=u_{_{\small{N}}}$ is indeed the unique global minimizer to the energy $E.$
\end{enumerate}

\begin{re}(i) To summarize the foregoing analysis, we have shown that for the full ball $B\ss\R^2$ and smooth boundary conditions, albeit with a topological change, there is a uniformly convex functional, which depends discontinuously on $x,$ but smoothly on $\grad u,$ such that the corresponding energy is uniquely globally minimised by a map that is everywhere Lipschitz but not $C^1(B)$.\\ 
%The fact that such a simple counterexample exists shows the rigidity of the incompressible case. In other words, admissible maps satisfying the constraint seem extremely rare.  \\

(ii) One might also be interested in this counterexample on the scale of Sobolev spaces. With this in mind, note that for any $N\in\N\sm\{1\}$ and $b\in\R^2$ we have\footnote{One might compare this to recent results given in \cite{RZZ18,RZZ20}. They obtain higher integrability, in particular, in the same spaces as we do, for counterexamples in compressible elasticity. Clearly their result is extraordinarily difficult, since the counterexamples are produced by means of convex integration.}
 \[u=\frac{R}{\sqrt{N}}e_{NR}+b\in W^{2,q}(B,\R^2)\sm W^{2,2}(B,\R^2) \;\mb{for any}\;1\le q<2.\]

Moreover, one might like to compare our result with the high-order regularity result, given in \cite{BOP92}.
They showed that for the special case of the Dirichlet functional and $u\in W^{2,q}(B,\R^2)$ with $q>2$ being a stationary point satisfying $\det\grad u=1$ a.e., then $u\in C^\infty(B,\R^2).$
It is possible that a similar result could be established, for a fairly general non-autonomous $p-$growth functional with the necessary changes in $q$. Let's assume for a second that such a result is indeed possible. Intriguingly, this seems to leave a `gap' at $q=2.$\\

(iii) The singular set $\Sigma$ in our example is, of course, just the origin $\Sigma=\{0\}.$  It remains an open question whether there can be other incompressible variational problems, including in incompressible elasticity, where the stationary points/minimizers possess a richer $\Sigma.$
\end{re}

\bibliography{LiteraturePhDMD}
\bibliographystyle{plain}
\vspace{1cm}
\textsc{Acknowledgements:} The first author is appreciative to the Department of Mathematics at the University of Surrey and was funded by the Engineering \& Physical Sciences Research Council (EPRSC).

\end{document}